\documentclass[a4paper,10pt,leqno]{article}
\usepackage{amsfonts,amssymb,amsmath,amsthm}
\usepackage{fullpage}
\usepackage{epic,eepic}
\usepackage{mathrsfs}
\newtheorem{theorem}{Theorem}[section]

\newtheorem{lemma}[theorem]{Lemma}

\theoremstyle{remark}

\numberwithin{equation}{section}

\newcommand{\A}{\mathbf A}
\newcommand{\N}{{\mathbb N}}
\newcommand{\R}{{\mathbb R}}

\newcommand{\K}{{\mathcal K}}
\newcommand{\eps}{{\varepsilon}}
\newcommand{\vark}{{\varkappa}}

\newcommand\g{\gamma}

\renewcommand\d{\delta}
\renewcommand\l{\lambda}
\newcommand\n{\nabla}
\renewcommand\O{\Omega}

\def\var{\varphi}

\def\vark{\varkappa}

\def\part{\partial}

\def\pr{^{\prime}}

\def\XXint#1#2#3{{\setbox0=\hbox{$#1{#2#3}{\int}$}
\vcenter{\hbox{$#2#3$}}\kern-.5\wd0}}

\def\W{\rlap{$\buildrel \circ \over W$}\phantom{W}}

\begin{document}

\title{
Potential estimates for quasi-linear parabolic equations }

\author{
{\large Vitali Liskevich}\\
\small Department of Mathematics\\
\small Swansea University\\
\small Swansea SA2 8PP, UK\\
{\tt v.a.liskevich@swansea.ac.uk}\\
\and
{\large Igor I.\,Skrypnik}\\
\small Institute of Applied\\
\small Mathematics and Mechanics\\
\small Donetsk 83114,  Ukraine\\
{\tt iskrypnik@iamm.donbass.com}
\and
{\large Zeev Sobol}\\
\small Department of Mathematics\\
\small  Swansea University\\
\small Swansea SA2 8PP, UK\\
{\tt z.sobol@swansea.ac.uk}
}

\date{}

\maketitle

\setlength{\unitlength}{0.0004in}
\begingroup\makeatletter\ifx\SetFigFont\undefined%
\gdef\SetFigFont#1#2#3#4#5{%
  \reset@font\fontsize{#1}{#2pt}%
  \fontfamily{#3}\fontseries{#4}\fontshape{#5}%
  \selectfont}%
\fi\endgroup%
\renewcommand{\dashlinestretch}{30}

\begin{abstract}
For a class of divergence type quasi-linear degenerate parabolic  equations with a Radon measure on the right hand side
we derive pointwise estimates for solutions via  nonlinear Wolff potentials.
\end{abstract}

\bigskip

\section{Introduction and main results}
In this note we give a parabolic extension of a by now classical result by Kilpel\"ainen-Mal\'y estimates \cite{KiMa}
who proved
pointwise estimates of solutions to quasi-linear $p$-Laplace type elliptic equations with measure in the right hand side. The estimates are
expressed in terms of the nonlinear Wolff potential of the right hand side. These estimates were subsequently extended to fully nonlinear
and subelliptic quasi-linear equation by Trudinger and Wang \cite{TW}. For the parabolic equations the corresponding result was recently given in \cite{DM2,DM3}, but
 only
for the "linear" case $p=2$ . Here we provide the estimates for parabolic equations in the degenerate case $p>2$.

Let $\O$ be a domain in $\R^n$, $T>0$. Let $\mu$ be a Radon measure on $\O$.
We are concerned with pointwise estimates for a class of non-homogeneous divergence type quasi-linear parabolic equations
of the type
\begin{equation}
\label{e0}
u_t- {\rm div}\,\A(x,t,u,\n u)=\mu\quad \text{in}\ \O_T=\O\times (0,T), \quad \O\subset \R^n,
\end{equation}
and assume that the following structure conditions are satisfied:
\begin{eqnarray}
\label{e1.3}
\A (x,t,u, \zeta)\zeta &\ge& c_1 |\zeta|^p,\quad \zeta\in \R^n,\\
\label{e1.4}
|\A (x,t,u, \zeta)|&\le& c_2 \zeta|^{p-1},
\end{eqnarray}
with some positive constants $c_1, c_2$,
%
whose model involves the parabolic $p$-Laplace equation
\begin{equation}
\label{e1.1} u_t-\Delta_pu=\mu ,\quad (x,t)\in
\O_T.
\end{equation}

\bigskip

Before formulating the main results, let us remind the reader of the definition of a weak solution
to equation \eqref{e1.1}.

We say that $u$ is a weak solution to \eqref{e1.1} if $u\in
V(\O_T):=C([0,T]; L^2_{loc}(\O))\cap L^p_{loc}(0,T;
W^{1,p}_{loc}(\O))$ and for any compact subset $\K$ of $\O$ and
any interval $[t_1,t_2]\subset (0,T)$ the integral identity
\begin{equation}
\label{e1.12b} \int_{\K}u\var dx\Big|_{t_1}^{t_2}+\int_{t_1}^{t_2}
\int_\K \left\{-u\var_\tau+\A(x,t,u,\n u) \n\var
\right\}dx\,d\tau=\int_{t_1}^{t_2}\int_{\K} \var \mu(dx)\,d\tau.
\end{equation}
for any $\var \in W^{1,2}_{loc}(0,T;L^2(\K))\cap L^p_{loc}(0,T; \W^{1,p}(\K))$.

Further on, we assume that $u_t\in L^2_{loc}(\O_T)$, since otherwise
we can pass to Steklov averages.

\medskip
The crucial role in our results is played by the truncated version
of the Wolff potential defined by
\begin{equation}
\label{wolff}
W^\mu_p (x,R)=\int_0^R \left(\frac{\mu(B_r(x))}{r^{n-p}}\right)^\frac1{p-1}\frac{dr}r.
\end{equation}

In the sequel, $\g$ stands for a constant which depends only on $n,p,c_1,c_2$ which may vary from line to line.

The main result of this paper  is the following theorem.

\begin{theorem}
\label{thm1.1b}Let $u$ be a weak solution to equation~\eqref{e0}.
For every $\l\in (0, \frac1n]$ there exists $\g>0$ depending on $n,c_1,c_2$ and $\l$, such that
for almost all $(y,s)\in \O_T$
and for $\rho\in (0,1)$ such that $B_{2\rho}(y)\times (s-4\rho^2,s+4\rho^2)\subset \O_T$ one has
\[
u(y,s)\le
\g\left\{\left(\frac1{\rho^{p+n}}\iint_{B_\rho\times(s-\rho^p,s+\rho^p)} u_+^{(1+\l)(p-1)}dxdt\right)^\frac{1}{1+\l(p-1)}
+1+ W^\mu_p(y, 2\rho)
\right\}.
\]
\end{theorem}

The estimate above is not homogeneous in $u$ which is usual for
such type of equations \cite{DiB, DiGV}. The proof of
Theorem~\ref{thm1.1b} is based on a suitable modifications of
De Giorgi's iteration technique \cite{DG} following the
adaptation of Kilpel\"ainen-Mal\'y technique \cite{KiMa} to
parabolic equations with ideas from \cite{LS2, Skr1}.


\bigskip




The rest of the paper contains the proof of the theorem.
\section{ Proof of Theorem~\ref{thm1.1b}}
\label{bdd}

%

\bigskip

We start with some auxiliary integral estimates for the solutions of \eqref{e0} which are formulated in the next lemma.

Define
\[
G(u)=
\left\{
\begin{array}{lll}
u&\text{for}\ &\ u>1,\\
u^{2-2\l}&\text{for}\ &\ 0<u\le 1.
\end{array}
\right.
\]
Set
\[
Q_\rho^{(\d)}(y,s)=B_\rho(y)\times
(s-\d^{2-p}\rho^p,\, s+\d^{2-p}\rho^p)\subset \O_T, \quad\rho\le R.
\]
\begin{lemma}
\label{lem2.2b}
Let the conditions of Theorem~\ref{thm1.1b} be fulfilled. Let $u$ be a solution to \eqref{e1.1}.
Then there exists a constant $\g>0$ depending only on $n,p,c_1,c_2$ such that for any $\eps\in(0,1),\, l,\d>0$,
any cylinder $Q_\rho^{(\d)}(y,s)$
and any  $\xi\in C_0^\infty(Q_\rho^{(\d)}(y,s))$ such that $\xi(x,t)=1$ for $(x,t)\in Q_{\rho/2}^{(\d)}(y,s)$
\begin{eqnarray}
\nonumber
&&\d^2 \int_{L(t)}G\left(\frac{u(x,t)-l}{\d}\right)\xi(x,t)^kdx
+\iint_{L}
\left(1+\frac{u-l}{\d}\right)^{-1+\l}\left(\frac{u-l}{\d}\right)^{-2\l}|\n u|^p\xi(x,\tau)^kdx\,d\tau\\
\nonumber
&\le& \g \d^2 \iint_L  \left(\frac{u-l}{\d}\right)|\xi_t|\xi^{k-1}dxd\tau
+\g \frac{\d^p}{\rho^p}\iint_L \left[\left(1+\frac{u-l}{\d}\right)^{1-\l}\left(\frac{u-l}{\d}\right)^{2\l}\right]^{p-1} \xi^{k-p}dxd\tau
\\
&&\ \ \ \ \ \ \ \ \ \ \ \ \  \  +\g\,  \d^{3-p}\mu(B_\rho(y))
, \label{e2.3b}
\end{eqnarray}
where $L=Q_\rho^{(\d)}(y,s)\cap\{u>l\}$, $L(t)=L\cap \{\tau=t\}$ and $\l\in (0,1)$,
$k>p$.
\end{lemma}
\proof First, note that
\begin{equation}
\label{e2.4b}
\int_l^u\left(1+\frac{s-l}{\d}\right)^{-1+\l}\left(\frac{s-l}{\d}\right)^{-2\l}ds\le \g \d,
\end{equation}
and
\begin{eqnarray}
\nonumber
\int_l^u dw\int_l^w\left(1+\frac{s-l}{\d}\right)^{-1+\l}\left(\frac{s-l}{\d}\right)^{-2\l}ds
\nonumber
=\int_l^u\left(1+\frac{s-l}{\d}\right)^{-1+\l}\left(\frac{s-l}{\d}\right)^{-2\l}(u-s)ds\\
\nonumber
\ge\frac12 (u-l)\int_l^{\frac{u+l}2}\left(1+\frac{s-l}{\d}\right)^{-1+\l}\left(\frac{s-l}{\d}\right)^{-2\l}ds
=\frac{\d^2}2\left(\frac{u-l}{\d}\right)\int_0^{\frac{u-l}{2\d}}(1+z)^{-1+\l}z^{-2\l}dz\\
\label{e2.5b}
\ge \g\d^2G\left(\frac{u-l}{\d}\right).
\end{eqnarray}
Test \eqref{e1.12b} by $\var$ defined by
\begin{equation}
\var(x,t)=\left[\int_l^{u(x,t)}\left(1+\frac{s-l}{\d}\right)^{-1+\l}\left(\frac{s-l}{\d}\right)^{-2\l}ds\right]_+\xi(x,t)^k,
\end{equation}
and $t_1=s-\d^{2-p}\rho^p$, $t_2=t$.
Using the Young inequality and \eqref{e2.4b} we have for any $t>0$
\begin{eqnarray*}
&&\int_{L(t)}\int_l^u dw\int_l^w\left(1+\frac{s-l}{\d}\right)^{-1+\l}\left(\frac{s-l}{\d}\right)^{-2\l}ds \xi^k dx\\
&+&\iint_L
\left(1+\frac{u-l}{\d}\right)^{-1+\l}\left(\frac{u-l}{\d}\right)^{-2\l}|\n u|^p\xi^k dxdt\\
&\le& \g \iint_L \int_l^udw \int_l^w\left(1+\frac{s-l}{\d}\right)^{-1+\l}\left(\frac{s-l}{\d}\right)^{-2\l}ds|\xi_t|\xi^{k-1}dxdt\\
&+&\g \iint_L \left[\left(1+\frac{u-l}{\d}\right)^{1-\l}\left(\frac{u-l}{\d}\right)^{2\l}\right]^{p-1} \xi^{k-p}dx\,dt
+\g \d^{3-p}\rho^p \mu(B_\rho(y)).
\end{eqnarray*}
From this using \eqref{e2.4b} and \eqref{e2.5b} we obtain the required~\eqref{e2.3b}.~\qed

\bigskip
Now set
\begin{equation}
\label{e2.6b}
\psi(x,t)=\frac1{\d}\left[\int_l^{u(x,t)}\left(1+\frac{s-l}{\d}\right)^{-\frac{1-\l}{p}}\left(\frac{s-l}{\d}\right)^{-\frac{2\l}{p}}ds\right]_+.
\end{equation}
The next lemma is a direct consequence of Lemma~\ref{lem2.2b}.
\begin{lemma}
\label{lem2.3b}
Let the conditions of Lemma~\ref{lem2.2b} be fulfilled.
Then
\begin{eqnarray}
\nonumber
&&\int_{L(t)}G\left(\frac{u-l}{\d}\right)\xi^k
dx+\d^{p-2}\iint_L |\n\psi|^p \xi^k dxd\tau\\
&&\le
\g \frac{\d^{p-2}}{\rho^p}\iint_L \left(1+\frac{u-l}{\d}\right)^{(1-\l)(p-1)}
\left(\frac{u-l}{\d}\right)^{2\l(p-1)}\xi^{k-p}dxd\tau\ \
\label{e2.10b}
+\g \frac{\rho^p}{\d^{p-1}}\mu({B_\rho(y)}). \ \
\end{eqnarray}
\end{lemma}

\bigskip

Let $(y,s)$ be an arbitrary point in $\O_T$.
Let
$
R\le \frac12\min\left\{1,{\rm dist}\,(y,\partial\O),
s^\frac12, (T-s)^\frac12\right\}
$
and
$
Q_R(y,s)=B_R(y)\times (s-R^2,s+R^2).
$
Fix $\rho\le R$ and
for $j=0,1,2,\dots$ set
\[
\rho_j=\rho 2^{-j},\quad Q_j=B_j\times (s-\d_j^{2-p}\rho_j^p,s+\d_j^{2-p}\rho_j^p),
\quad B_j=B_{\rho_j}(y),\quad L_j=Q_j\cap\O_T\cap\{u(x,t)>l_j\}.
\]
Let $\xi_j\in C_0^\infty(Q_j)$ be such that $\xi_j(x,t)=1$ for $(x,t)\in B_{j+1}\times (s-\frac34 \d_j^{2-p}\rho_j^p,\,s+\frac34 \d_j^{2-p}\rho_j^p)$,
$|\n \xi_j|\le \g \rho_j^{-1}$, $|\frac{\partial \xi_j}{\partial t}|\le \g \d_j^{p-2}\rho_j^{-p}$.

The sequences of positive numbers $(l_j)_{j\in\N}$ and $(\d_j)_{j\in\N}$ are defined inductively as follows.

Set $l_0=0$ and assume that $l_1,l_2,\dots,l_j$ and $\d_0,\d_1,\dots,\d_{j-1}$ have been already chosen in such a way that
$\d_k=l_{k+1}-l_k$.
Let us show  how to chose $l_{j+1}$ and $\d_j$.


For 
$l\ge l_j+\rho_j$
set
\begin{eqnarray}
\nonumber
A_j(l)=\frac{(l-l_j)^{p-2}}{\rho_j^{n+p}}\iint_{\widetilde L_j}\left(\frac{u-l_j}{l-l_j}\right)^{(1+\l)(p-1)}\xi_j^{k-p}dxd\tau\, \\
\label{e3.1b}
+\sup_{|t-s|\le (l-l_j)^{2-p}\rho_j^p}\, \frac1{\rho_j^n}\int_{\widetilde L_j(t)}G\left(\frac{u-l_j}{l-l_j}\right)\xi_j^k dx,
\end{eqnarray}
where $\widetilde L_j=\widetilde Q_j\cap \O_T\cap\{u(x,t)>l_j\}$,
$\widetilde Q_j=B_j\times (s-(l-l_j)^{2-p}\rho_j^p,s+(l-l_j)^{2-p}\rho_j^p)$.


Fix a number $\varkappa\in (0,1)$ depending on $n,p,c_1,c_2$, which will be specified later. Set $\hat\d_0=\max\{1,\rho_0\}$, $\hat \d_j=\rho_j$.
For $j=0,1,2,\dots$, if
\begin{equation}
\label{e3.3b}
A_j(l_j+\hat\d_j)\le \varkappa,
\end{equation}
we set $\d_j=\hat\d_j$ and $l_{j+1}=l_j+\d_j$.

Note that $A_j(l)$ is continuous as a function of $l$ and
$A_j(l)\searrow 0$ as $l\to \infty$. So if
\begin{equation}
\label{e3.4b}
A_j(l_j+\hat\d_j)> \varkappa,
\end{equation}
there exists $\bar l>l_j+\hat\d_j$ such that $A_j(\bar l)=\varkappa$. In this case we set
$l_{j+1}=\bar l$ and $\d_j=l_{j+1}-l_j$.

Note that our choices guarantee that $\widetilde Q_j\subset Q_R(y,s)$ and
\begin{equation}
\label{e3.5b}
A_j(l_{j+1})\le \varkappa.
\end{equation}
The following lemma is a key in the Kilpel\"ainen-Mal\'y technique \cite{KiMa}.
\begin{lemma}
\label{lem3.1b}
Let the conditions of Theorem~\ref{thm1.1b} be fulfilled.
There exists $\g>0$ depending on the data, such that for all $j\ge 1$
we have
\begin{equation}
\label{e3.6b}
\d_j\le \frac12 \d_{j-1}+\rho_j
+\g \left(\frac1{\rho_j^{n-p}}\mu(B_j)
\right)^\frac1{p-1}.
\end{equation}
\end{lemma}
\proof Fix $j\ge 1$. Without loss assume that
\begin{equation}
\label{e3.7b}
\d_j>\frac12 \d_{j-1},\quad \d_j>\rho_j,
\end{equation}
since otherwise \eqref{e3.6b} is evident. The second inequality in \eqref{e3.7b} guarantees that $A_j(l_{j+1})=\varkappa$ and $\widetilde Q_j=Q_j$.

Next we claim that under conditions \eqref{e3.7b} there is a $\g>0$ such that
\begin{equation}
\label{claim}
\d_j^{p-2} \rho_j^{-(p+n)}|L_j|\le \g \vark.
\end{equation}
Indeed,
for $(x,t)\in L_j$ one has
\begin{equation}
\label{Lj}
\frac{u(x,t)-l_{j-1}}{\d_{j-1}}=1+\frac{u(x,t)-l_{j}}{\d_{j-1}}\ge 1.
\end{equation}
Note  that the first inequality in \eqref{e3.7b} yields
$\xi_{j-1}=1$ on $Q_j$. Hence
\begin{eqnarray*}
&&\d_j^{p-2} \rho_j^{-(p+n)}|L_j|\le \d_j^{p-2} \rho_j^{-(p+n)}\iint_{L_j} G\left(\frac{u-l_{j-1}}{\d_{j-1}}\right)\xi_{j-1}^{k}dx\,d\tau\\
&\le& \rho_j^{-n}\sup_{|t-s|\le \d^{2-p}_j \rho_j^p}\int_{L_j(t)}G\left(\frac{u-l_{j-1}}{\d_{j-1}}\right)\xi_{j-1}^kdx
\le 2^n \rho_{j-1}^{-n}\sup_{|t-s|\le \d^{2-p}_{j-1} \rho_{j-1}^p}\int_{L_{j-1}(t)}G\left(\frac{u-l_{j-1}}{\d_{j-1}}\right)\xi_{j-1}^kdx\le 2^n\vark,
\end{eqnarray*}
which proves the claim.

Let us estimate the terms in the right hand side of \eqref{e3.1b} with $l=l_{j+1}$. For this we decompose $L_j$ as
$L_j=L_j\pr\cup L_j^{\prime\prime}$,
\begin{equation}
\label{decomp}
L_j\pr=\left\{(x,t)\in L_j\,:\,\frac{u(x,t)-l_j}{\d_j}<\eps \right\}, \quad L_j^{\prime\prime}=L_j\setminus L\pr_j,
\end{equation}
where $\eps\in(0,1)$ depending on $n,p,c_1,c_2$ is small enough to be determined later.

 By \eqref{claim} we have
\begin{equation}
\label{e3.9b}
\frac{\d_j^{p-2}}{\rho_j^{n+p}}\iint_{L\pr_j}\left(\frac{u-l_j}{\d_j}\right)^{(1+\l)(p-1)}\xi_j^{k-p}dxd\tau
\le \frac{\d_j^{p-2}}{\rho_j^{n+p}}\iint_{L\pr_j}\eps^{(1+\l)(p-1)}dx\,d\tau
\le 2^{n} \eps^{(1+\l)(p-1)} \varkappa.
\end{equation}

Set
\begin{equation}
\label{psi-j}
\psi_j(x,t)=\frac1{\d_j}\left(\int_{l_j}^{u(x,t)}\left(1+\frac{s-l_j}{\d_j}\right)^{-\frac{1-\l}{p}}
\left(\frac{s-l_j}{\d_j}\right)^{-\frac{2\l}{p}}ds\right)_+,
\end{equation}
and
\[
\rho(\l)=\frac{p}{p-1-\l}.
\]

Note that
$\l\le \frac1n$ due to the assumption.

The following inequalities  are easy to verify
\begin{eqnarray}
\label{psi}
c\psi_j(x,t)^{\rho(\l)}\le \left(\frac{u(x,t)-l_j}{\d_j}\right) \ \text{for}\ (x,t)\in L_j, \ \text{and} \\
\left(\frac{u(x,t)-l_j}{\d_j}\right)
\le c(\eps)\psi_j(x,t)^{\rho(\l)},
\quad (x,t)\in L_j^{\prime\prime}.
\end{eqnarray}
%

Hence
\begin{eqnarray}
\nonumber
\frac{\d_j^{p-2}}{\rho_j^{n+p}}\iint_{L_j^{\prime\prime}}\left(\frac{u-l_j}{\d_j}\right)^{(1+\l)(p-1)}\xi_j^{k-p}dxd\tau\\
\label{e3.10b}
\le
\g(\eps)\frac{\d_j^{p-2}}{\rho_j^{n+p}}\iint_{L_j^{\prime\prime}}  \psi_j^{p\frac{n+\rho(\l)}{n}}\xi_j^{k-p}dxd\tau.
\end{eqnarray}
The integral in the second terms of the right hand side of \eqref{e3.10b} is estimated by using
the Gagliardo--Nirenberg inequality in the form \cite[Chapter~II,Theorem~2.1]{LaSU}
as follows
\begin{eqnarray}
\nonumber
&&\g(\eps)\frac{\d_j^{p-2}}{\rho_j^{n+p}}\iint_{L_j^{\prime\prime}}  \psi_j^{p\frac{n+\rho(\l)}{n}}\xi_j^{k-p}dxd\tau\\
\label{e3.11a}
&&\le \g\,\left(\sup_{|t-s|\le \d_j^{2-p}\rho_j^p}\frac1{\rho_j^n}\int_{L_j(t)}\psi_j^{\rho(\l)}dx\right)^\frac{p}n
\left(\frac1{\rho_j^n}\iint_{L_j}\left|\n \left(\psi_j \xi_j^\frac{(k-p)n}{p(n+\rho(\l))}\right)\right|^p dx\,d\tau\right).
\end{eqnarray}

Let us estimate separately the first factor in the right hand side of \eqref{e3.11a}.
\begin{eqnarray}
\nonumber
&&\sup_{|t-s|\le \d_j^{2-p}\rho_j^p}\int_{L_j(t)}\psi_j^{\rho(\l)}dx\ \stackrel{\text{by \eqref{psi}}}
{\le} c^{-1}\sup_{|t-s|\le \d_j^{2-p}\rho_j^p}\int_{L_j(t)}\frac{u-l_j}{\d_j}\,dx\\
\nonumber
&&\stackrel{\text{by \eqref{e3.7b}}}
{\le} 2c^{-1}\sup_{|t-s|\le \d_j^{2-p}\rho_j^p}\int_{L_j(t)}\frac{u-l_{j-1}}{\d_{j-1}}\xi_{j-1}^k\,dx\\
\nonumber
&&\stackrel{\text{by \eqref{Lj}}}
\le 2c^{-1}\sup_{|t-s|\le \d_{j-1}^{2-p}\rho_{j-1}^p}\int_{L_{j-1}(t)}G\left(\frac{u-l_{j-1}}{\d_{j-1}}\right)\xi_{j-1}^k\,dx\\
&&\stackrel{\text{by \eqref{e3.5b}}}
{\le} 2c^{-1}\vark \rho_{j-1}^n = \g \rho_j^n \vark.
\label{psi1}
\end{eqnarray}

Combining \eqref{e3.10b}, \eqref{e3.11a} and \eqref{psi1} we obtain
\begin{eqnarray}
\nonumber
\frac{\d_j^{p-2}}{\rho_j^{n+p}}\iint_{L_j^{\prime\prime}}\left(\frac{u-l_j}{\d_j}\right)^{(1+\l)(p-1)}\xi_j^{k-p}dxd\tau\\
\le
\g(\eps)\vark^\frac{p}n \, \d_j^{p-2}\rho_j^{-n}\iint_{L_j}\left|\n \left(\psi_j \xi_j^\frac{(k-p)n}{p(n+\rho(\l))}\right)\right|^p dx\,d\tau.
\end{eqnarray}
For the last term in the above inequality we estimate by \eqref{claim} and \eqref{psi}
\begin{eqnarray}
\nonumber
&&\d_j^{p-2}\rho_j^{-n}\iint_{L_j}\psi_j^p\left|\n  \xi_j\right|^p dx\,d\tau
\le \g \d_j^{p-2}\rho_j^{-n-p}\iint_{L_j}\psi_j^p dx\,d\tau\\
\nonumber
&& \le \g  \d_j^{p-2}\rho_j^{-n-p}\iint_{L_j} \left(\frac{u-l_j}{\d_j}\right)^{p-1-\l}dx\,d\tau
\stackrel{\text{by \eqref{e3.7b}}}{\le}\g  \d_{j-1}^{p-2}\rho_j^{-n-p}\iint_{L_j} \left(\frac{u-l_{j-1}}{\d_{j-1}}\right)^{p-1-\l}\xi_{j-1}^{k-p}dx\,d\tau\\
\label{psi-p}
&&\le \g  \d_{j-1}^{p-2}\rho_j^{-n-p}\iint_{L_{j-1}} \left(\frac{u-l_{j-1}}{\d_{j-1}}\right)^{(p-1)(1+\l)}\xi_{j-1}^{k-p}dx\,d\tau
\le \g\vark.
\end{eqnarray}

By Lemma~\ref{lem2.3b}
\begin{eqnarray}
\nonumber
\frac1{\rho^n_j}\int_{L_j(t)}G\left(\frac{u-l_j}{\d_j}\right)\xi_j^k
dx+\frac{\d_j^{p-2}}{\rho_j^n}\iint_{L_j} |\n\psi_j|^p \xi_j^k dxd\tau\\
\nonumber
\le
\g \frac{\d_j^{p-2}}{\rho_j^{p+n}}\iint_{L_j} \left(1+\frac{u-l_j}{\d_j}\right)^{(1-\l)(p-1)}
\left(\frac{u-l_j}{\d_j}\right)^{2\l(p-1)}\xi_j^{k-p}dxd\tau\ \ &&\\
\label{e3.d}
+\g \frac{\rho_j^{p-n}}{\d_j^{p-1}}\mu({B_{\rho_j}(y)}). \ \ &&
\end{eqnarray}

 Using the decomposition \eqref{decomp} and the first inequality in \eqref{e3.7b} we have
\begin{eqnarray}
\nonumber
&&
\d_j^{p-2}\rho_j^{-(n+p)}\iint_{L_j}\left(1+\frac{u-l_j}{\d_j}\right)^{(1-\l)(p-1)}\left(\frac{u-l_j}{\d_j}\right)^{2\l(p-1)}dx\,d\tau
\\
\nonumber
&& \le \g \eps^{2\l(p-1)}\d_j^{p-2}\rho_j^{-(n+p)}|L_j|+\g(\eps)\d_{j-1}^{p-2}\rho_{j-1}^{-(n+p)}\iint_{L_{j-1}}\left(\frac{u-l_{j-1}}{\d_{j-1}}\right)^{(1+\l)(p-1)}dx\,d\tau
 \\
 \label{A}
 && \ \ \ \ \ \ \ \ \le \g \eps^{2\l(p-1)} \vark +\g(\eps) \vark.
\end{eqnarray}
Thus we obtain the following estimate for the first term of $A_j(l_{j+1})$:
\begin{eqnarray}
\nonumber
&&\frac{\d_j^{p-2}}{\rho_j^{n+p}}\iint_{L_j}\left(\frac{u-l_j}{\d_j}\right)^{(1+\l)(p-1)}dxd\tau\\
\label{first}
&&\le \g
\eps^{2\l(p-1)}\vark +\g(\eps)\vark^\frac pn\left(\vark +\d_j^{1-p}\rho_j^{p-n}\mu(B_j)\right).
\end{eqnarray}

Let us estimate the second term in the right hand side of \eqref{e3.1b}.
By
\eqref{e3.d} we have
\begin{eqnarray}
\nonumber
&&\sup_{|t-s|\le \d_j^{2-p}\rho_j^p}\rho_j^{-n}\int_{L_j(t)} G\left(\frac{u-l_j}{\d_j}\right)\xi_j^k dx\\
\nonumber
&&\le \d_j^{p-2}\rho_j^{-(n+p)}\iint_{L_j}\left(1+\frac{u-l_j}{\d_j}\right)^{(1-\l)(p-1)}\left(\frac{u-l_j}{\d_j}\right)^{2\l(p-1)}\xi_j^{k-p}dx\,d\tau
+ \g \d_j^{1-p}\rho_j^{p-n}\mu(B_j)\\
\nonumber
&&\text{(by using the decomposition \eqref{decomp} and \eqref{first})}\\
\label{second}
&&\le \g\eps^{2\l(p-1)}\vark +\g(\eps )\vark^\frac pn\left(\vark +\d_j^{1-p}\rho_j^{p-n}\mu(B_j)\right)
+ \g \d_j^{1-p}\rho_j^{p-n}\mu(B_j).
\end{eqnarray}

Combining \eqref{A} and \eqref{second} and choosing $\eps$ appropriately we can find $\g_1$ and $\g$ such that
\begin{equation}
\label{kappa}
\vark \le \g_1\vark^\frac pn\left(\vark +\d_j^{1-p}\rho_j^{p-n}\mu(B_j)\right) + \g \d_j^{1-p}\rho_j^{p-n}\mu(B_j).
\end{equation}
Now choosing $\vark<1$ such that
$\displaystyle \vark^\frac pn =\frac1{2\g_1}$ we have
\begin{equation}
\label{deltaj}
\d_j\le \g \left(\rho_j^{p-n}\mu(B_j)\right)^\frac1{p-1},
\end{equation}
which completes the proof of the lemma.\qed
%

In order to complete the proof of Theorem~\ref{thm1.1b} we sum up \eqref{e3.6b} with respect to $j$ from 1 to $J-1$
\begin{eqnarray}
\nonumber
&&l_J\le \g \d_0 +\g \sum_{j=1}^\infty \rho_j + \g \sum_{j=1}^\infty \left(\rho_j^{p-n}\mu(B_j)\right)^\frac1{p-1}\\
\label{lJ}
&&\le \g (\d_0+\rho +W^\mu_p(y, 2\rho)).
\end{eqnarray}

Let us estimate $\d_0$. There are two cases to consider.
If $l_1=\hat\d_0=\max\{1,\rho\}$ then
$\d_0=\max\{1,\rho\}$.
If on the other hand $l_1$ and $\d_0$ are defined by $A_0(l_1)=\varkappa$ then by \eqref{e3.1b}
\[
\vark=\frac{\d_0^{p-2}}{\rho^{n+p}}\iint_{Q_\rho^{(\d_0)}}\left(\frac{u_+}{\d_0}\right)^{(1+\l)(p-1)}\xi_0^{k-p} dx\,d\tau
+\sup_{|t-s|<\d_0^{2-p}\rho^p}\rho^{-n}\int_{B_\rho}G\left(\frac{u_+}{\d_0}\right)\xi_0^k\,dx.
\]
Using the decomposition \eqref{decomp} with $\eps$ chosen via $\vark$,
and Lemma~\ref{lem2.3b} one can see that
\[
\sup_{|t-s|<\d_0^{2-p}\rho^p}\rho^{-n}\int_{B_\rho}G\left(\frac{u_+}{\d_0}\right)dx \le \vark/2+
\frac{\d_0^{p-2}}{\rho^{n+p}}\iint_{Q_\rho^{(\d_0)}}\left(\frac{u_+}{\d_0}\right)^{(1+\l)(p-1)}dx\,d\tau.
\]
Note that $\d_0\ge \max\{1,\rho\}$, thus $\d_0^{2-p}\rho^p\le \rho^p$.
Hence we obtain
\[
\vark\le \g \frac{\d_0^{p-2}}{\rho^{n+p}}\iint_{Q_\rho^{(\rho)}}\left(\frac{u_+}{\d_0}\right)^{(1+\l)(p-1)} dx\,d\tau.
\]
Combining this with the first case we have
\begin{equation}
\label{e3.17b}
\d_0\le \g  \left\{\left(\frac1{\rho^{p+n}}\iint_{B_\rho(y)\times(s-\rho^p,s+\rho^p)} u_+^{(1+\l)(p-1)}dxdt\right)^\frac{1}{1+\l(p-1)}
+1+\rho
\right\}.
\end{equation}

Hence the sequence $(l_j)_{j\in\N}$ is convergent, and $\d_j\to 0\,(j\to\infty)$, and we can pass to the limit $J\to\infty$ in \eqref{lJ}.
Let $l=\lim_{j\to\infty}l_j$. From \eqref{e3.5b} we conclude that
\begin{equation}
\label{e3.19b}
\frac1{\rho_j^{n+p}}\iint_{Q_j}(u-l)_+^{(1+\l)(p-1)}\le \g \vark\,  \d_j^{1+\l(p-1)} \to 0\quad (j\to \infty).
\end{equation}
Choosing $(y,s)$ as a Lebesgue point of the function $(u-l)_+^{(1+\l)(p-1)}$ we conclude that $u(y,s)=l$ and hence $u(y,s)$
is estimated from above by

\[
u(y,s)\le
\g\left\{\left(\frac1{\rho^{p+n}}\iint_{B_\rho(y)\times(s-\rho^p,s+\rho^p)} u_+^{(1+\l)(p-1)}dxdt\right)^\frac{1}{1+\l(p-1)}
+1+\rho
+ W^\mu_p(y, 2\rho)
\right\}
\]
Applicability of the Lebesgue differentiation theorem
follows from \cite[Chap.~II,\,Sec.~3]{Guz}.



\section*{Acknowledgment}
The authors would like to thank Giuseppe Mingione for useful discussion and for providing
a preprint of \cite{DM3} prior to publication.

%
%
%
%
%
%
%
%
%

\begin{small}

\end{small}


\end{document}